\documentclass[reqno,12pt]{amsart}
\usepackage{xr}
\usepackage{graphicx}
\usepackage[usenames,dvipsnames,svgnames,table]{xcolor}
\usepackage{amsrefs}

\usepackage{epstopdf}
\usepackage{booktabs}
\usepackage{amssymb,amsmath}
\usepackage{amsthm}
\usepackage{enumerate}
\newtheorem{theorem}{Theorem}[section]

\newtheorem{definition}[theorem]{Definition}

\numberwithin{equation}{section}

\newcommand{\A}{{\mathcal{A}}}
\newcommand{\CA}{{\mathcal{A}}}
\newcommand{\CB}{{\mathcal{B}}}
\newcommand{\CL}{{\mathcal{L}}}
\newcommand{\CO}{{\mathcal{O}}}

\renewcommand{\ll}{{\langle}}

\newcommand{\rr}{{\rangle}}
\newcommand{\g}{\mathfrak{g}}
\newcommand{\C}{{\mathbb C}}
\newcommand{\Z}{{\mathbb Z}}

\newcommand{\Cone}{{\mathrm{Cone}}}
\newcommand{\rep}{{\mathrm{rep}}}
 \newcommand{\geo}{{\mathrm{geo}}}
 
 \newcommand{\Tr}{{\mathrm{Tr}}}
\newcommand{\End}{{\mathrm{End}}}
\newcommand{\Todd}{{\mathrm{Todd}}}

 \renewcommand{\c}{{\mathfrak{c}}}
\begin{document}
\title{
The equivariant Riemann-Roch theorem and the graded Todd class}

\author{ Mich{\`e}le Vergne}
\address{Universit\'e Denis-Diderot-Paris 7, Institut de Math\'ematiques de Jussieu,
 C.P.~7012\\ 2~place Jussieu,   F-75251 Paris~cedex~05}
\email{michele.vergne@imj-prg.fr}

\maketitle

\begin{abstract}
Let $G$ be a torus with Lie algebra $\g$ and
 let $M$ be a $G$-Hamiltonian  manifold  with Kostant line bundle $\CL$ and proper moment map. Let $\Lambda\subset \g^*$ be the weight lattice  of $G$.
 We consider a parameter $k\geq 1$ and the multiplicity  $m(\lambda,k)$
 of the quantized representation $RR_G(M,\CL^k)$ .
Define   $\ll\Theta(k),f\rr=
\sum_{\lambda\in \Lambda} m(\lambda,k) f(\lambda/k)$
for $f$ a test function on $\g^*$.
We prove that the distribution $\Theta(k)$ has an asymptotic development
$\ll\Theta(k),f\rr\sim k^{\dim M/2}\sum_{k=0}^{\infty} k^{-n} \ll DH_n,f\rr$
where the distributions $DH_n$ are the twisted Duistermaat-Heckman distributions associated to the graded Todd class of $M$.
When $M$ is compact, and $f$ polynomial, the asymptotic series is finite and exact.

\end{abstract}

\section{Introduction}

Let $G$ be a torus with Lie algebra $\g$.
Identify $\hat G$ to a lattice $\Lambda$ of  $\g^*$.
If $\lambda\in \Lambda$, we denote by  $g^{\lambda}$
the corresponding character of $G$. If $g=\exp (X)$ with $X\in \g$, then
$g^{\lambda}=e^{i\ll\lambda,X\rr}$.

Let $M$  be a prequantizable $G$-hamiltonian manifold with symplectic form $\Omega$,
 Kostant line bundle $\CL$,  and moment map $\Phi:M\to \g^*$.
Assume $M$  compact and of dimension $2d$.
The  Riemann-Roch quantization  $RR_G(M,\CL)$ is a
virtual finite dimensional representation of $G$, constructed as the index of a Dolbeaut-Dirac operator on $M$.
 The dimension of
the space $RR_G(M,\CL)$ will be called the Riemann-Roch number of $(M,\CL)$.
 The character of the representation  of $RR_G(M,\CL)$ is a function on $G$, denoted by  $RR_G(M,\CL)(g).$
 We write $$RR_G(M,\CL)(g)=\sum_{\lambda \in\Lambda} m_{\rep}(\lambda)g^{\lambda}.$$

{\em The typical example is the case where $M$ is a projective  manifold, and
 $\CL$ the corresponding ample bundle. Then $$RR_G(M,\CL)(g)=\sum_{i=0}^{d} (-1)^i \Tr_{H^i(M,\CO(\CL))}(g)$$ is the alternate sum of the traces of the action of $g$ in
the cohomology spaces of
 $\CL$.
 In particular $\dim RR_G(M,\CL)=\sum_{i=0}^{d} (-1)^i \dim H^i(M,\CO(\CL))$
  is given by the Riemann-Roch formula.}

\bigskip

It is natural to introduce the $k^{th}$ power $\CL^k$ of the line bundle $\CL$.
  Thus $$RR_G(M,\CL^k)(g)=\sum_{\lambda \in\Lambda} m_{\rm rep}(\lambda,k)g^{\lambda}.$$

Assume $k\geq 1$.
We associate to $(M,\CL)$  the distribution on $\g^*$ given by
$$\ll\Theta_M(k), f\rr=\sum_{\lambda\in \Lambda} m_{\rm rep}(\lambda,k) f(\lambda/k),$$
where $f$ is  a test function.

\bigskip
{\em Example: When $M$ is a toric manifold associated to the Delzant polytope $P$,
 then $\dim RR_G(M,\CL)$ is the number of integral points in  the convex polytope $P$,
 and $\frac{1}{k^d}\ll\Theta_M(k), f\rr$
 is the Riemann sum of the values of $f$ on the sample points $\frac{\Lambda}{k}\cap P$.}

\bigskip

We prove that
$\Theta_M(k)$ has an asymptotic behavior when the integer  $k$ tends to $\infty$
of the form
$$\Theta_M(k)\sim k^d\sum_{n=0}^{\infty} k^{-n} DH_n$$
where $DH_n$ are distributions on $\g^*$ supported on $\Phi(M)$.
We determine the  distributions $DH_n$  in terms of
the decomposition of the equivariant Todd class  ${\Todd}(M)$ of $M$
in its homogeneous components  ${\Todd}_n(M)$
in the graded equivariant cohomology ring of $M$.
The distribution $DH_0$ is the Duistermaat-Heckmann measure.
The asymptotics are exact when $f$ is a polynomial.
This generalizes the weighted Ehrhart polynomial for an integral polytope,
and the asymptotic behavior  of Riemann sums  over convex integral polytopes
established by Guillemin-Sternberg \cite{guisteriemann}.

\bigskip

We then consider the case where $M$  is a prequantizable $G$-hamiltonian manifold, not necessarily compact, but
with proper moment map $\Phi: M\to \g^*$.
The formal quantization of $(M,\CL^k)$ (\cite{weitsman})  is defined by
$$RR_G(M,\CL^k)(g)=\sum_{\lambda\in \Lambda} m_{\rm geo}(\lambda,k)g^{\lambda}.$$
Here $m_{\rm geo}(\lambda,k)$  is the geometric multiplicity function constructed by
Guillemin-Sternberg
in terms of the Riemann-Roch number  of the reduced fiber
  $M_\lambda=\Phi^{-1}(\lambda)/G$ of the moment map.
When $M$ is compact,  Meinrenken-Sjamaar \cite{MeiSja} proved that $m_{\rep}(\lambda,k)=m_{\geo}(\lambda,k)$,
so this purely geometric definition extends the definition of
$RR_G(M,\CL^k)$ given in terms of  index theory  when $M$ is compact.

Similarly, we construct distributions $DH_n$ on $\g^*$ using
 the equivariant cohomology classes ${\Todd}_n(M)$ and push-forwards by the proper map $\Phi$.
The main result of this announcement is that the distribution
$\Theta_M(k)$  defined by
$$\ll\Theta_M(k), f\rr=\sum_{\lambda\in \Lambda} m_{\rm geo}(\lambda,k) f(\lambda/k),$$
is asymptotic to
$k^d\sum_{n=0}^{\infty} k^{-n} DH_n.$

A similar result holds for  Dirac operators twisted by powers of a line bundle $\CL^k$.

Recall that we introduced a truncated Todd class (of the cotangent bundle $T^*M$)  for determining the multiplicities of the equivariant index of any transversally elliptic operator on $M$ \cite{vergneserre}.
Here the use of the parameter $k$ allows us to have families of such equivariant indices, and the full series $\sum_{n=0}^{\infty} \Todd_n(M)$  enters in
 the description of the asymptotic behavior. This is similar to the Euler-Mac Laurin formula evaluating sums of the values of a function
 at integral points of an interval involving all Bernoulli numbers.
We finally give some information on the piecewise polynomial behavior of the distributions $DH_n$.

\section{Equivariant cohomology}
Let $N$ be a $G$-manifold and let $\A(N)$ be the space of differential forms on $N$,
graded by its exterior degree.
Following \cite{berver82} and \cite{witten82},
an equivariant form is a $G$-invariant smooth function
$\alpha: \g\to \A(N),$ thus
$\alpha(X)$ is a differential form on $N$ depending differentiably of $X\in \g$.
Consider the operator
\begin{equation}\label{DX}
d_\g\alpha(X)=d\alpha(X)-\iota(v_X)\alpha(X)
\end{equation}
 where $\iota(v_X)$ is the contraction by the vector field $v_X$
generated by the action of $-X$ on $N$.
Then $d_\g$ is an odd operator with square $0$,
and  the equivariant cohomology is defined
to be the cohomology space of $d_\g$.
It is important to note that the dependance of $\alpha$ on $X$ may be $C^\infty$.
If the dependance of $\alpha$ in $X$ is polynomial, we denote by
$H^*_G(N)$ the corresponding  $\Z$-graded algebra.
By definition, the grading of $P(X)\otimes \mu$, $P$ a homogeneous polynomial and $\mu$ a differential form on $N$,
is the exterior degree of $\mu$ plus twice the polynomial degree in $X$.

The Hamiltonian structure on $M$ determines the equivariant symplectic form
$\Omega(X)=\ll \Phi,X\rr + \Omega$.

Choose a  $G$-invariant Riemannian metric on $M$.
This  provides the  tangent bundle $TM$ with the structure of a Hermitian vector bundle.
Let $J(A)=\det_{\C^d} \frac{e^{A}-1}{A}$, an invariant  function of  $A\in \End(\C^d)$.
 Then, $J(0)=1$. Consider $\frac{1}{J(A)}$ and its Taylor expansion at $0$:
 $$\frac{1}{J(A)}=\det_{\C^d}(\frac{A}{e^{A}-1})= \sum_{n=0}^{\infty}  B_{n}(A).$$
 Each  function $B_{n}(A)$ is an invariant polynomial of degree $n$ on $\End(\C^d)$ and by the Chern Weil construction,
  $B_{n}$
determines an equivariant characteristic class
${\Todd}_{n}(M)(X)$ on $M$ of homogeneous degree $2n$. Remark that ${\Todd}_0(M)=1$.
We define the formal series of equivariant cohomology classes:
   $${\Todd}(M)(X)=\sum_{n=0}^{\infty} {\Todd}_{n}(M)(X).$$
   For $X$ small enough,
   the series is convergent, and ${\Todd}(M)(X)$  is the equivariant Todd class  of  $M$.
   In particular ${\Todd}(M)(0)$ is the usual Todd class of $M$.

In the rest of this note, using the Lebesgue measure $d\xi$ determined by the lattice $\Lambda$, we may identify distributions and generalized functions on $\g^*$
and we may write $\ll \theta,f\rr=\int_{\g^*}\theta(\xi)f(\xi)d\xi$ for the value of a distribution $\theta$ on a test function $f$ on $\g^*$.

\section{ The compact case}
Let $M$ be a compact $G$-Hamiltonian manifold.
Recall (see \cite{ber-get-ver}) the "delocalized Riemann-Roch formula". For $X\in \g$ sufficiently small, we have
$$RR_G(M,\CL)(\exp X)=\frac{1}{(2i\pi)^{d}}\int_{M} e^{i\Omega(X)}\Todd(M)(X).$$
Here $i=\sqrt{-1}$.

For each integer $n$, consider
the analytic function on $\g$ given by
$$\theta_n(X)=\frac{1}{(2i\pi)^{d}}\int_{M} e^{i \Omega(X)}{\Todd}_n(M)(X).$$
Thus when $X\in \g$ is small, then $\sum_{n=0}^{\infty}\theta_n(X)$ is a
convergent series with sum the equivariant Riemann-Roch index $RR_G(M,\CL)(\exp X)$.
When $n=0$, $$\theta_0(X)=\frac{1}{(2i\pi)^{d}}\int_{M} e^{i \Omega(X)}$$
is the equivariant volume of $M$, and the Fourier transform $DH_0$  of
$\theta_0$ is the Duistermaat-Heckmann measure of $M$,
a piecewise polynomial distribution on $\g^*$.

\begin{theorem}\label{mytheo1}
Let $DH_n$ be the Fourier transform of $\theta_n$.
Then $DH_n$ is a  distribution supported on $\Phi(M)$.
 For any polynomial function $P$ of degree $N$ on $\g^*$, we have
$$\sum_{\lambda\in \Lambda}m_{{\rm rep}}(\lambda)
P(\lambda)=\sum_{n\leq N+d} \int_{\g^*}DH_n(\xi)P(\xi)d\xi.$$

In particular, we have  the following   Euler-MacLaurin formula for the Riemann-Roch number of $(M,\CL)$:
$$\dim RR_G(M,\CL)=
\sum_{\lambda\in \Lambda}m_{{\rm rep}}(\lambda)=\int_{\g^*} \sum_{n\leq d} DH_n(\xi)d\xi.$$
\end{theorem}

We now give a theorem for smooth functions.

\begin{theorem}\label{mytheo2}
When the integral parameter $k$ tends to  $\infty$,
 the distribution $\Theta_M(k)$ admits the asymptotic expansion
 $$\Theta_M(k)\sim k^{d}\sum_{n=0}^{\infty} k^{-n} DH_n.$$
\end{theorem}
Let us sketch the proof of Theorems \ref{mytheo1} and \ref{mytheo2}.
It is easy to see
that the distributions $DH_n$ are supported on the image $\Phi(M)$ of $M$ by the moment map.
Furthermore, it follows from the piecewise quasi-polynomial behavior of the function
$m_{\rm rep}(\lambda,k)$ that for $P$ a homogeneous polynomial of degree $N$,
the sum  $\sum_{\lambda\in \Lambda}m_{{\rm rep}}(\lambda,k)
P(\lambda)$ is a quasi-polynomial function of $k\geq 1$
of degree less or equal than $N+d$. Thus
 Theorem \ref{mytheo1} will be a consequence of Theorem \ref{mytheo2} that we now prove.

The Fourier transform of
$\Theta_M(k)$ is
$$\sum_{\lambda\in \Lambda} m_{\rm rep}(\lambda,k)e^{i\ll \lambda,X/k\rr}=
 RR_G(M,\CL^k)(\exp (X/k)).$$
 Against a test function $\phi$ of $X$,
 when $k$ is large, this is
$$\frac{1}{(2i\pi)^{d}}\int_{M\times \g} e^{i k\Omega+ik\Phi(X/k)}\Todd(M)(X/k)\phi(X) dX$$
 $$=\sum_{n,m}
\frac{1}{(2i\pi)^{d}}\int_{M} e^{i\Phi(X)} \frac{1}{m!}k^{m} (i\Omega)^{m}
 \Todd_{n}(M)(X/k)\phi(X) dX.$$
 For each $m$, only the term of differential degree $2d-2m$ of ${\Todd}_n(M)$ contributes to the integral, and this term
  is homogeneous in $X$ of degree $n+m-d$. This implies the result.

\bigskip

Let us give an example of the asymptotic expansion.

Let $P_1(\C)$ equipped with the torus action $g([z_1,z_2])=[gz_1,z_2]$ of $g=e^{i\theta}$, in homogeneous coordinates.
 We consider $M=P_1(\C)\times P_1(\C)$ with diagonal action, and let $\CL$ be its Kostant line bundle.
 Then we have
$$RR(M,\CL^k)(g)=\sum_{j\in \Z} m_{\rm rep}(j,k) g^j $$
with
$$
m_{\rm rep}(j,k)=\begin{cases}
0\hspace{30mm} {\rm if}\,  j<-2k ,\\
    2k+1+j\hspace{14mm} {\rm if}\ -2k\leq j\leq 0,\\
    2k+1-j\hspace{14mm} {\rm if}\ 0\leq j\leq 2k,\\
    0\hspace{30mm} {\rm if} \,  j>2k .\\
    \end{cases}
$$

We have

$$\Theta(k)\sim k^2 (DH_0+\frac{1}{k} DH_1+ \frac{1}{k^2}DH_2+ \frac{1}{k^3}DH_3+\cdots.)$$
Let us give the explicit formulae for $DH_0,DH_1,DH_2,DH_3$.

 $$\ll DH_0,f\rr=\int_{-2}^2 m(\xi) f(\xi) d\xi$$
with
$$
m(\xi)=
\begin{cases}
   2+\xi\hspace{14mm} {\rm if}-2\leq \xi\leq 0,\\
    2-\xi\hspace{14mm} {\rm if}\ 0\leq \xi\leq 2,\\
\end{cases}
$$

$$\ll DH_1,f\rr=\int_{-2}^2 f(\xi) d\xi,$$

$$\ll DH_2,f\rr=\frac{5}{12}f(-2)+\frac{1}{6}f(0)+\frac{5}{12}f(2),$$

$$\ll DH_3,f \rr=-\frac{1}{12}f'(-2)+ \frac{1}{12}f'(2).$$

\bigskip

We now sketch another proof of Theorem \ref{mytheo2} which
can be extended to the non compact case. We use Paradan's decomposition
(\cite{parlocfirst},\cite{parloc}, see also \cite{sze-ver})
of $RR_G(M,\CL)$ in a sum of simpler  characters supported on cones.
Let us consider a  generic value $r$ of the moment map,
and choose a scalar product on $\g^*$.
Then there exists a certain finite subset $\CB(r)$ of   $\g^*$, and for each
$\beta\in \CB(r)$,
a cone $C(\beta)$ in $\g^*$
and an (infinite dimensional) representation $P_{\beta,k}$
such that
$$RR_G(M,\CL^k)=\sum_{\beta\in \CB(r)} P_{\beta,k}.$$
 Here $P_{\beta,k}(g)=\sum_{\lambda\in \Lambda \cap k C(\beta)} m_{{\rm rep},\beta}(\lambda,k) g^\lambda$.
Thus $\Theta(k)$ is decomposed in $\sum_{\beta\in \CB(r)} \Theta_\beta(k).$
Similarly each distribution $DH_n$
is decomposed as
$DH_n=\sum_{\beta\in \CB(r)} DH_{n,\beta}$
and the support of $DH_{n,\beta}$ is contained in the cone $C_\beta$.
It is easily verified that, for each $\beta$,
the distribution $\Theta_\beta(k)$ is asymptotic to
$k^d\sum_{n=0}^{\infty} k^{-n} DH_{n,\beta}$.
Here we use  the explicit Euler Mac Laurin expansion on half lines,
and  convolutions of such distributions. The proof is entirely similar to
the case of a polytope given in \cite{ber-ver-elm}.

Let us return to the example of the case of $M=P_1(\C)\times P_1(\C)$, for $r<0$ a small
negative number. Then
$\CB(r)=\{-2,r,0,2\}$.
We have
$$P_{\beta=-2,k}(g)=-\sum_{j<-2k} (2k+1+j) g^j,$$
$$P_{\beta=r,k}(g)=\sum_{j=-\infty}^{j=\infty}   (2k+1+j) g^j,$$
$$P_{\beta=0,k}(g)=-2\sum_{j>0} j g^j,$$
$$P_{\beta=2,k}(g)=\sum_{j>2k}   (j-(2k+1)) g^j.$$

Consider, for example, the asymptotic development of
the distribution
$$\ll\Theta_{\beta=2}(k),f\rr=\sum_{j>2k} (j-(2k+1)) f(j/k).$$
It is easy to see that this distribution is the convolution $K(k)*K(k)$
where $K(k)$ is the distribution defined by $\ll K(k),f\rr=\sum_{j>k} f(j/k)=\sum_{j\geq k} f(j/k)-f(1).$
We then use the explicit exact Euler Mac-Laurin formula to evaluate  the distribution $K(k)$, thus its convolution.
In particular the Fourier transform of $K(k)*K(k)$ coincides with
the analytic function $\frac{e^{2ix}}{(1-e^{-ix/k})^2}$ for $(1-e^{-ix/k})\neq 0$.
As is  natural, the asymptotic series of distributions $q^{-d}\sum_{n=0}^{\infty}q^{n}DH_{\beta=2,n}$
is the unique series of distributions supported on $\xi\geq 2$ and with Fourier transform, for $x\neq 0$,
the Laurent series in $q$ of $\frac{e^{2ix}}{(1-e^{-iqx})^2}$ at $q=0$.

\section{Proper moment maps}

Consider the case where $M$  is non necessarily compact, but  $\Phi: M\to \g^*$ is a proper map.
One can then define (\cite{weitsman},\cite{parformal1}) the formal geometric
quantification of $M$  with respect to the line bundle $\CL^k$
 to be
 $$RR_{G,{\rm geo}}(M,\CL^k)(g)=\sum_{\lambda\in \Lambda} m_{\rm geo}(\lambda,k) g^\lambda,$$
using a function $m_{\rm geo}(\xi)$ on $\g^*$ .
The  definition of the function  $m_{\rm geo}(\xi)$  is due to Guillemin-Sternberg \cite{Gui-ste}.
Let us recall its delicate definition (\cite{MeiSja}, see also \cite{pep-vergnespin1}). There is a closed set $\CA$,  union of affine hyperplanes,  such that
if $r$ is in the complement of $\CA$, then either $r$ is not in $\Phi(M)$
or $r$ is a regular value of $\Phi$.
Consider the open subset $\g^*_{{\rm reg}}=\g^*\setminus \CA$.
 When $\xi\in \g^*_{{\rm reg}}$ but not in $\Phi(M)$, $m_{\rm geo}(\xi)$ is defined to be $0$.
 If  $\xi\in \g^*_{{\rm reg}}\cap \Phi(M)$, the reduced fiber $M_\xi=\Phi^{-1}(\xi)/G$  is a compact symplectic orbifold, and $m_{\rm geo}(\xi)$
  is defined to be a sum of  integrals on the various strata of the compact orbifold $M_\xi$.
When $\lambda\in \g^*_{{\rm reg}}\cap \Lambda\cap \Phi(M)$, then $M_\lambda$ is a prequantizable  compact symplectic orbifold and $m_{\rm geo}(\lambda)$
 is the Riemann-Roch number of $M_\lambda$ equipped with its Kostant orbifold line bundle.
Let $\lambda\in \Lambda$ be any point in  $\Phi(M)$.
Choose a vector $\epsilon$, such that
$\lambda+t\epsilon$ is in  $\Phi(M)\cap \g^*_{{\rm reg}}$ for any $t>0$ and sufficiently small.
It  can be proved, using the wall crossing formulae of Paradan \cite{parwallcrossing},
 that $(\lim_\epsilon m_{\rm geo})(\lambda)=\lim_{t>0, t\to 0}m_{\rm geo}(\lambda+t\epsilon)$
 is independent of the choice of  such an $\epsilon$.
 This allows us to define  $m_{\rm geo}(\lambda)$ by "continuity on $\Phi(M)$" for any $\lambda\in \Lambda$.

The $[Q,R]=0$ theorem  (\cite{MeiSja}, \cite{ma-zhang}, \cite{parformal2}) asserts
 that $RR_{G,{\rm geo}}(M,\CL)$ coincides with a representation of $G$ defined  using
index theory. In particular  $RR_{G,{\rm geo}}(M,\CL)$ coincides with
 $RR_G(M,\CL)$ when $M$ is compact.
 However in the rest of this note, we only use the geometric definition of
$RR_{G,{\rm geo}}(M,\CL)$.

Replacing $\CL$ by $\CL^k$, and the moment map $\Phi$ by $k\Phi$,
define the distribution, with parameter $k$,
$$\ll\Theta_M(k), f\rr=\sum_{\lambda\in \Lambda} m_{\rm geo}(\lambda,k) f(\lambda/k).$$

As in the compact case, the asymptotic behavior of $\Theta_M(k)$ is determined by the graded Todd class, using push-forwards by the proper map $\Phi$.
Indeed if $\alpha$  is an equivariant cohomology class with polynomial coefficients, then
the Duistermaat-Heckman twisted distribution $DH(M,\Phi,\alpha)$ is well defined by
the formula
$$\ll DH(M,\Phi,\alpha),f\rr=\frac{1}{(2i\pi)^d}
\int_{M\times \g} e^{i\Omega(X)} \alpha(X) \hat f(X) dX$$
where $\hat f(X)=\int_{\g^*} e^{i\ll \xi,X\rr}f(\xi)d\xi$ is the Fourier transform of the test function $f(\xi)$ (see  \cite{con-pro-ver-inf}).
It is a distribution supported on $\Phi(M)$.

\begin{definition}
We define $DH_n$ to be the distribution on $\g^*$ associated to the equivariant cohomology class $\Todd_n(M)$:
$$\ll DH_n,f\rr=\frac{1}{(2i\pi)^d}
\int_{M\times \g} e^{i\Omega(X)} \Todd_n(X) \hat f(X) dX.$$
\end{definition}

The distribution $DH_0$ is the Duistermaat-Heckman measure, a locally polynomial function.

The distribution $DH_n$ is given by a polynomial function
on  each connected component of the  open set $\g^*_{{\rm reg}}$. Its restriction to
$\g^*_{{\rm reg}}$ vanishes when $n> d-\dim G$.
Furthermore, if all stabilizers of points of $M$ are connected,
it follows from Witten non abelian localization theorem
that $$m_{\rm geo}(\lambda,k)=k^d\sum_{n=0}^{\infty} k^{-n} DH_n(\lambda/k)$$
when $\lambda/k$ is a regular value of $\Phi$.  Otherwise, it can be defined by limit of the function
$m_{\rm geo}(\xi,k)=k^d\sum_{n=0}^{\infty} k^{-n} DH_n(\xi/k)$, along
 $\xi=\lambda+t\epsilon_\lambda$ and $t>0, t\to 0$,
where the direction $\epsilon_\lambda$ is chosen to be arbitrary
if $\lambda$ does not belong to $\Phi(M)$,
or in such a way that $\lambda+t\epsilon_\lambda$ stays in $\Phi(M)$ if $\lambda\in \Phi(M)$. Similar formulae can be given without assumption on connected stabilizers.

We can see that,  for any $n$, the distributions $DH_n$
can be expressed (but not uniquely)
as derivatives of  locally polynomial
functions associated to symplectic submanifolds $M^{T}$
 where $T$ are  subtori of $G$  .

The main result of this note is the following theorem.
\begin{theorem}
When the integer $k$ tends to $\infty$,
 $$\Theta_M(k)\sim k^d\sum_{n=0}^{\infty} k^{-n} DH_n.$$
\end{theorem}

Let us  sketch  the proof of this theorem, in the case where each stabilizer is connected.
We use Paradan decomposition formula \cite{parloc}, \cite{parlocfirst}.
We choose $r$ a generic element of $\g^*_{{\rm reg}}$.
As in the compact case, there is a locally finite set
$\CB(r)\subset \Phi(M)$, cones $C_\beta$, and  decompositions
$$DH_n=\sum_{\beta\in \CB(r)} DH_{n,\beta}$$
where $DH_{n,\beta}$ are supported on $C_\beta$.
The functions $DH_{n,\beta}$ are  given by polynomial functions on each connected component of  $\g^*_{{\rm reg}}$
and vanishes on $\g^*_{\rm reg}$ when $n>d-\dim G$.
Thus the locally polynomial function $A_\beta(\xi,k)=
k^d\sum_{n=0}^{\infty} k^{-n} DH_{n,\beta}(\xi/k)$ is well defined when $\xi/k\in \g^*_{\rm reg}$.
For each $\beta\in \CB(r)$,
choose a direction $\epsilon_\beta$ such that
$\beta+t\epsilon_\beta$ is in $\Phi(M)\cap \g^*_{{\rm reg}}$ for $t>0$ small.
Then $w_{\beta}(\lambda,k)=\lim_{t>0, t\to 0} A_\beta(\lambda+t \epsilon_\beta,k)$
is well defined.
Define $$P_{\beta, k}(g)=\sum_{\lambda\in \Lambda} w_\beta(\lambda,k) g^{\lambda}$$
and $$\ll\Theta_{\beta,{\rm geo}}(k),f\rr=\sum_{\lambda\in \Lambda} w_\beta(\lambda,k) f(\lambda/k).$$
As before,  it is easy to see that
$\Theta_{\beta,{\rm geo}}(k)\sim  k^d \sum_{n=0}^{\infty} k^{-n} DH_{n,\beta}$.
Here we use the following "continuity" result on partition function (see for example
\cite{con-pro-ver-box}).  Let $\Delta$ be a unimodular  list of non zero vectors in $\Lambda$, and $\gamma\in \g$ generic.
   There is a unique function $K$ (the Kostant partition function)  on $\Lambda$ supported
   on the half space $\ll\xi,\gamma\rr\geq 0$ and such that
   $\sum_{\lambda\in \Lambda} K(\lambda) g^{\lambda}=
   \prod_{\alpha\in \Delta}\frac{1}{1-g^{\alpha}}$ for $g$ in the open set
   $\prod_{\alpha\in \Delta}(1-g^{\alpha})\neq 0$.
Let $d=|\Delta|$.
   Consider the Laurent series expansion  in $q$
  $$ \prod_{\alpha\in \Delta}\frac{1}{1-e^{q\ll\alpha,X\rr}}=
  q^{-d}\sum_{n=0}^{\infty} q^n U_n(X)$$
  and the distributions $D_n$ on $\g^*$ supported on
   the half space $\ll\xi,\gamma\rr\geq 0$, such that
   $$\int_{\g^*} D_n(\xi) e^{i\ll\xi,X\rr}=U_n(X)$$
   when $\prod_\alpha \ll\alpha,X\rr\neq 0$.
Define $T(\xi)=\sum_{n=0}^\infty
D_n(\xi)$, which is well defined outside a system of hyperplanes.
Then for any $\lambda\in \Lambda$,
and $\epsilon_\Delta$  generic and belonging to the cone $\Cone(\Delta)$ generated by $\Delta$,
we have $K(\lambda)=\lim_{t>0, t\to 0} T(\lambda+t\epsilon_\Delta)$.

Define $P_{r,k}=\sum_{\beta\in \CB(r)} P_{\beta,k}$.
It remains to see that $P_{r,k}=RR_{G,{\rm geo}}(M,\CL^k)$.
This is not immediate, since
we do not have a global representation theoretic object for describing
$RR_{G,{\rm geo}}(M,\CL^k)$.
Each coefficient
$m_{\rm geo}(\lambda,k)$ is defined using a limit direction depending on $\lambda$
while   each $w_\beta(\lambda,k)$ is defined using the same limit direction (depending on $\beta$) for any $\lambda$. So additivity is not clear.
However, we can prove that $P_{r,k}$ is independent of $r$, using \cite{parwallcrossing}.
This is  very similar to the technique used in \cite{agapito} to establish decompositions  \`a la Paradan of characteristic functions of polyhedra.
It then follows that $P_{r,k}=RR_{G,{\rm geo}}(M,\CL^k)$. Indeed for each connected component
$\c$
of $\g^*_{{\rm reg}}$ contained in $\Phi(M)$, we choose $r$ in $\c$.
In the decomposition $P_{r,k}=\sum_{\beta\in \CB(r)} P_{\beta,k}$, the
  term $w_\beta(\lambda,k)$  for $\beta=r\in \CB(r)$
is  the  polynomial function
coinciding with $m_{\rm geo}(\lambda,k)$ for $\lambda\in k \overline{\c}$.
The other terms $w_\beta(\lambda,k)$ for $\beta\in \CB(r)$ and $\beta\neq r$ vanishes when $\lambda\in k\overline{\c}$ (\cite{parloc}, see also \cite{sze-ver}).

A quicker route, but less instructive, for determining asymptotics of
$\Theta_{M,{\rm geo}}$ would be to take a test function with small support around a point $r\in \g^*$. Then we can choose $\epsilon_\lambda$ coinciding with $\epsilon_\beta$ for all $\beta\in \CB(r)$ and in the support of the test function $f$. The additivity is  immediate on those $\beta$.

%To illustrate the  method and the delicate direction problem, we look  at a very simple two dimensional case. Let $\R^2$ with standard basis $e_1,e_2$.
%Let $S=\{(\xi_1,\xi_2), -1\leq \xi_1\leq 1, \xi_2\geq 0\}$.
%Consider $M=P_1(\C)\times \C$
%with action of the two dimensional torus
%$(g_1,g_2)([z_1,z_2],z)= ([g_1z_1,z_2],g_2z_2)$
%and moment map such that $\Phi(M)=S$. Here $g_1=e^{i\theta_1}, g_2=e^{i\theta_2}$.
%The function $m_{\rm geo}(\xi,1)$ is the characteristic function of
%the set $S\cap \Z^2$.
%Choose  $r=t e_2$ on the  vertical line, and $t<0$.
% Then the elements $ \beta\in \CB(r)$ are
% $\beta_0=r$, $\beta_1=[0,0]$, $\beta_2=[-1,0]$, $\beta_3=[0,1]$.
%The corresponding directions can be chosen as
%$\epsilon_0=e_2$, $\epsilon_1=e_2$, $\epsilon_2=e_2+e_1$,$\epsilon_3=-e_1+e_2$
%and there is no common choice.
%
%
%
%
%
%  The corresponding functions are given by
% $$w_{\beta_0,k=1}(\xi)=0,$$
%$$w_{\beta_1}(\lambda_1,\lambda_2)=1$ if $\lambda_2\geq 0$
%
%$$w^1_{\beta_2}(\lambda_1,\lambda_2)=-1$$ if $\lambda_2\geq 0, \lambda_1<-1$
%
%$w^1_{\beta_2}(\lambda_1,lambda_2)=-1$ if $\lambda_2\geq 0, \lambda_1>1$.
%
%
%
%Thus the sum is indeed equal to $m_{\rm geo}(\lambda)$.
%However, $m_{\rm geo}([-2,0])=0$ while both $w^1_{\beta_0}$ and $w^1_{\beta_1}$ do not vanish.
%
%
%
%
%

\end{document}